\documentclass[10pt]{article}
\usepackage{amsmath,amssymb,amsthm,amsfonts,enumerate,times,epsfig,color,wrapfig,graphicx,tikz,subfig,pdfpages,mathrsfs,hyperref,cleveref}
\usepackage[numbers]{natbib}
\usepackage[letterpaper,twoside,outer=1.3in,vmargin=1.3in,]{geometry}
\usepackage{hyperref}
\hypersetup{
    colorlinks,
    citecolor=black,
    filecolor=black,
    linkcolor=black,
    urlcolor=black
}

\newtheorem{theorem}{Theorem}[section]
\newtheorem{CO}[theorem]{Corollary}
\newtheorem{LE}[theorem]{Lemma}

\newtheorem{DE}[theorem]{Definition}
\newcounter{claim_nb}[theorem]
\newtheorem*{claim*}{Claim}

\newcommand{\ignore}[1]{}
\newcommand{\conv}{\mathrm{conv}}

\newcommand{\core}{\mathrm{core}}

\newcommand{\setcore}[1]{\mathrm{setcore}\!\left(#1\right)}

\title{The rainbow covering number of clean tangled clutters}

\author{Ahmad Abdi \and G\'{e}rard Cornu\'{e}jols}

\begin{document}
	
\maketitle
	
\begin{abstract}
In this brief note, we prove a min-min equality for a clean tangled clutter, that the rainbow covering number is equal to the connectivity of its setcore.\\

\end{abstract}

\section{Introduction}\label{sec:intro}

A clutter $\mathcal{C}$ over ground set $V$ is \emph{tangled}, or \emph{$2$-cover-minimal}, if it has covering number $2$, and every element belongs to a minimum cover~\cite{Abdi-kwise}. Denote by $G(\mathcal{C})$ the graph over vertex set $V$ whose edges correspond to the minimum covers of $\mathcal{C}$. A cover of $\mathcal{C}$ is \emph{rainbow} if it intersects every connected component of $G(\mathcal{C})$ at most once. The \emph{rainbow covering number of $\mathcal{C}$}, denoted $\mu({\mathcal{C}})$, is the minimum size of a rainbow cover of $\mathcal{C}$; if there is no rainbow cover, then $\mu({\mathcal{C}}):=+\infty$. 

A clutter is \emph{clean} if it has no minor that is a delta or the blocker of an extended odd hole. Let $\mathcal{C}$ be a clean tangled clutter over ground set $V$. It is known that $G(\mathcal{C})$ is a bipartite graph~\cite{Abdi-hole}. The \emph{core of $\mathcal{C}$}, denoted $\core(\mathcal{C})$, is the set of all members of $\mathcal{C}$ that intersect every minimum cover exactly once. It is known that the core of $\mathcal{C}$ corresponds uniquely to a set-system in $\{0,1\}^d$, where $d$ is the number of connected components of $G(\mathcal{C})$ -- this is defined in \S\ref{sec:prelim}. The set-system, denoted by $\setcore{\mathcal{C}}\subseteq \{0,1\}^d$, is called a \emph{setcore of $\mathcal{C}$}, and the convex hull of $\setcore{\mathcal{C}}$ is a full-dimensional polytope containing $\frac12 {\bf 1}$ in its interior~\cite{Abdi-progeo}.

Clean tangled clutters have been a subject of recent study. There is an intriguing interplay between the polyhedral geometry and the combinatorics of clean tangled clutters. For example, the convex hull of the setcore is a simplex if, and only if, the setcore is the cocycle space of a projective geometry over $GF(2)$~\cite{Abdi-progeo}. 

The \emph{connectivity} of a set-system $S\subseteq \{0,1\}^d$, denoted $\lambda(S)$, is the minimum number of variables that appear in a \emph{generalized set covering (GSC)} inequality valid for $\conv(S)$, i.e., an inequality of the form $\sum_{i\in I}x_i + \sum_{j\in J} (1-x_j)\geq 1$ for some $I,J\subseteq [d],I\cap J=\emptyset$; if $S=\{0,1\}^d$, then $\lambda(S):=+\infty$. 

In this brief note, we prove that the two parameters that we defined above are equal, further stressing the synergy between the polyhedral geometry and the combinatorics of clean tangled clutters:

\begin{theorem}\label{mono-connectivity}
For every clean tangled clutter $\mathcal{C}$, we have $\mu({\mathcal{C}}) = \lambda(\setcore{\mathcal{C}})$.
\end{theorem}

\section{Definitions and preliminaries}\label{sec:prelim}

\paragraph{Clutters} 

Let $V$ be a finite set of \emph{elements}, and let $\mathcal{C}$ be a family of subsets of $V$ called \emph{members} or \emph{sets}. $\mathcal{C}$ is a \emph{clutter} over \emph{ground set} $V$ if no member contains another~\cite{Edmonds70}. A \emph{cover} is a subset $B\subseteq V$ such that $B\cap C\neq \emptyset$ for all $C\in \mathcal{C}$. The \emph{covering number} of $\mathcal{C}$, denoted $\tau(\mathcal{C})$, is the minimum cardinality of a cover. A cover is \emph{minimal} if it does not contain another cover. The \emph{blocker} of $\mathcal{C}$, denoted $b(\mathcal{C})$, is the clutter over ground set $V$ whose members are the minimal covers of $\mathcal{C}$~\cite{Edmonds70}. It is well-known that $b(b(\mathcal{C}))=\mathcal{C}$~\cite{Isbell58,Edmonds70}. Take disjoint $I,J\subseteq V$. The \emph{minor} of $\mathcal{C}$ obtained after \emph{deleting $I$} and \emph{contracting~$J$}, denoted $\mathcal{C}\setminus I/J$, is the clutter over ground set $V-(I\cup J)$ whose members consist of the inclusion-wise minimal sets of $\{C-J:C\in \mathcal{C}, C\cap I=\emptyset\}$. 
It is well-known that $b(\mathcal{C}\setminus I/J) = b(\mathcal{C})/I\setminus J$~\cite{Seymour76}.

\paragraph{Clean clutters} 
Two clutters are \emph{isomorphic} if one is obtained from the other by relabeling the ground set. Take an integer $n\geq 3$. Denote by $\Delta_n$ the clutter over ground set $[n]:=\{1,\ldots,n\}$ whose members are $\{1,2\},\{1,3\},\ldots,\{1,n\}$, $\{2,3,\ldots,n\}$. Any clutter isomorphic to $\Delta_n$ is called a \emph{delta of dimension $n$}. A delta is equal to its blocker. 
Given an odd integer $n\geq 5$, an \emph{extended odd hole of dimension $n$} is any clutter whose ground set can be relabeled as $[n]$ so that its minimum cardinality members are precisely $\{1,2\},\{2,3\},\ldots,\{n-1,n\},\{n,1\}$. 
Recall that a clutter is clean if it has no minor that is a delta or the blocker of an extended odd hole. Testing cleanness of a clutter belongs to P~\cite{Abdi-int-rest}. 

\paragraph{Core and setcore} Let $\mathcal{C}$ be a clean tangled clutter over ground set $V$. Recall that $G:=G(\mathcal{C})$ is the graph over vertex set $V$ whose edges correspond to the minimum covers of $\mathcal{C}$. Recall that $G$ is a bipartite graph. Let $d$ be the number of connected components of $G(\mathcal{C})$, and for each $i\in [d]$, denote by $\{U_i,V_i\}$ the bipartition of the $i\textsuperscript{th}$ connected component of $G$. Recall that the core of $\mathcal{C}$, denoted $\core(\mathcal{C})$, is the set of all members of $\mathcal{C}$ that intersect every minimum cover exactly once. 

\begin{theorem}[\cite{Abdi-progeo}, Theorem 2.9]\label{core}
$\core(\mathcal{C})=\left\{C\in \mathcal{C}:C\cap (U_i\cup V_i)\in \{U_i,V_i\} \text{ for each } i\in [d]\right\}$. 
\end{theorem}

The \emph{setcore of $\mathcal{C}$ with respect to $(U_1,V_1;U_2,V_2;\ldots;U_d,V_d)$} is the set-system $S\subseteq \{0,1\}^d$ that has a point $p\in S$ for every $C\in \core(\mathcal{C})$ such that $p_i = 0$ if and only if $C\cap (U_i\cup V_i) = U_i$, for all $i\in [d]$. By \Cref{core}, the set-system $S$ is well-defined. We denote $S$ by $\setcore{\mathcal{C}:U_1,V_1;U_2,V_2;\ldots;U_d,V_d}$. As the reader can imagine, we will not use this notation often, and use $\setcore{\mathcal{C}}$ as short-hand notation. Note however that $\setcore{\mathcal{C}}$ is defined only up to isomorphism. 

\begin{theorem}[\cite{Abdi-progeo}, Theorem 1.5]\label{setcore}
$\conv(\setcore{\mathcal{C}})$ is a full-dimensional polytope contained in $[0,1]^d$ and containing $\frac12 {\bf 1}$ in its interior. In particular, $\setcore{\mathcal{C}}$ does not have duplicated coordinates, $\core(\mathcal{C})$ is nonempty, and has covering number two.
\end{theorem}

This immediately implies the following.

\begin{CO}\label{rank=1}
If $d\leq 2$, then $\setcore{\mathcal{C}}=\{0,1\}^d$.\qed
\end{CO}

We also need the following earlier result.

\begin{theorem}[\cite{Abdi-dyadic}, Theorem 2.5, and \cite{Abdi-progeo}, Lemma 2.6]\label{deletion-contraction}
Suppose $G$ is not a connected graph. Let $\{U,U'\}$ be the bipartition of a connected component of $G$. Then $\mathcal{C}\setminus U/U'$ is a clean tangled clutter such that $\core(\mathcal{C}\setminus U/U')\subseteq \core(\mathcal{C})\setminus U/U'$.
\end{theorem}

\section{Proof of \Cref{mono-connectivity}}\label{sec:rainbow-covering}

Let $\mathcal{C}$ be a clean tangled clutter over ground set $V$, let $G:=G(\mathcal{C})$, and let $d$ be the number of connected components of $G$. For each $i\in [d]$, let $\{U_i,V_i\}$ be the bipartition of the $i\textsuperscript{th}$ connected component of~$G$. Let $\mu:=\mu(\mathcal{C})$ be the rainbow covering number, and let $\lambda:=\lambda(\setcore{C})$ be the connectivity of $\setcore{C}$. We need a key notion.

\begin{DE}
A \emph{monochromatic cover} of $\mathcal{C}$ is a cover that is monochromatic in some proper $2$-vertex-coloring of $G$. A monochromatic cover of $\mathcal{C}$, say of the form $\bigcup_{i\in I} V_i$ for some $I\subseteq [d]$, is \emph{irreducible} if for each $j\in I$, $\left(\bigcup_{i\in I, i\neq j} V_i\right)\cup U_j$ is not a cover. 
\end{DE}

Observe that every rainbow cover is also monochromatic. We need the following parameters: \begin{enumerate}
\item[$\mu_1(\mathcal{C})$:] the minimum size of a cover of $\core(\mathcal{C})$ that is monochromatic in some proper $2$-vertex-coloring of $G$. It can be readily seen that any such cover that is inclusion-wise minimal intersects every component of $G$ at most once, and $\mu\geq \mu_1(\mathcal{C}) = \lambda$.
\item[$\mu_2(\mathcal{C})$:] the minimum number of components of $G$ intersected by a monochromatic cover of $\mathcal{C}$. 
\item[$\mu_3(\mathcal{C})$:] the minimum number of components of $G$ intersected by an irreducible monochromatic cover of $\mathcal{C}$.
\end{enumerate} 

Let $\mu_i:=\mu_i(\mathcal{C})$ for $i=0,1,2$. To prove \Cref{mono-connectivity}, it remains to show that $\mu_1\geq \mu$. The first lemma below is a minor extension of [\cite{Abdi-progeo}, Lemma 2.8] and the second lemma is related to [\cite{Abdi-progeo}, Theorem~5.2].

\begin{LE}\label{setcore0}
Suppose for some $u,v\in V$, every member of $\core(\mathcal{C})$ containing $u$ also contains $v$. Then $u,v$ belong to the same part of the bipartition of a connected component of $G$.
\end{LE}
\begin{proof}
By \Cref{core}, it suffices to show that $u,v$ belong to the same connected component of $G$. Suppose otherwise. In particular, $G$ is not connected. Let $\{U,U'\}$ be the bipartition of the connected component containing $u$ where $u\in U'$. Then $\mathcal{C}\setminus U/U'$ is a clean tangled clutter such that $\core(\mathcal{C}\setminus U/U')\subseteq \core(\mathcal{C})\setminus U/U'$ by \Cref{deletion-contraction}.
Let $w$ be a neighbor of $u$ in $G$; so $w\in U$. Then $\{w,u\}$ is a cover of $\mathcal{C}$. As every member of $\core(\mathcal{C})$ containing $u$ also contains $v$, it follows that $\{w,v\}$ is a cover of $\core(\mathcal{C})$, implying in turn that $\core(\mathcal{C})\setminus U/U'$ has $\{v\}$ as a cover. However, $\core(\mathcal{C}\setminus U/U')\subseteq \core(\mathcal{C})\setminus U/U'$, so $\core(\mathcal{C}\setminus U/U')$ has a cover of cardinality $1$, a contradiction to \Cref{setcore}.
\end{proof}

\begin{LE}\label{irred-mono}
If $V_1\cup \cdots\cup V_k$ is an irreducible monochromatic cover for some integer $k\in [r]$, then there exists a monochromatic minimal cover $B$ such that $B\subseteq\bigcup_{i=1}^k V_i$ and $|B\cap V_i|=1$ for each $i\in [k]$.
\end{LE}
\begin{proof}
Out of all the monochromatic minimal covers of $\mathcal{C}$ contained in $\bigcup_{i=1}^k V_i$, pick one of minimum cardinality, call it $B$. As $\bigcup_{i=1}^k V_i$ is an irreducible monochromatic cover, it follows that $B\cap V_i\neq \emptyset, i\in [k]$. To finish the proof of the lemma, it suffices to show that $|B\cap V_1|=1$. Suppose for a contradiction that $|B\cap V_1|\geq 2$. Let $I:=B-V_1$, $J:= V-(U_1\cup V_1\cup I)$, and $\mathcal{C}':=\mathcal{C}\setminus I/J$, a minor over ground set $U_1\cup V_1$. Assume in the first case that $\tau(\mathcal{C}')\geq 2$. Then $\mathcal{C}'$ is clean and tangled, and $G[U_1\cup V_1]\subseteq G(\mathcal{C}')$. Thus $G(\mathcal{C}')$ is a connected bipartite graph whose bipartition is inevitably $\{U_1,V_1\}$. It therefore follows from \Cref{rank=1} that $U_1,V_1\in \mathcal{C}'$. However, $B\cap V_1 = B-I$ is a cover of $\mathcal{C}'$ disjoint from $U_1$, a contradiction. Assume in the remaining case that $\tau(\mathcal{C}')\leq 1$. That is, there is a $D\in b(\mathcal{C})$ such that $D\subseteq U_1\cup V_1\cup I$ and $|D-I|\leq 1$. As $D\subseteq (V_1\cup \cdots \cup V_k)\cup U_1$, and $V_1\cup \cdots \cup V_k$ is an irreducible monochromatic cover, it follows that $D\subseteq \bigcup_{i=1}^k V_i$. But then $D$ is a monochromatic minimal cover of $\mathcal{C}$ contained in $\bigcup_{i=1}^k V_i$ and $$|D| = |D-I| + |D\cap I| \leq 1 + |B-(U_1\cup V_1)|<|B\cap (U_1\cup V_1)| + |B-(U_1\cup V_1)| = |B|,$$ a contradiction to our minimal choice of $B$. As a result, $|B\cap V_1|=1$, as desired.
\end{proof}

\begin{LE}\label{mu-mu-mu}
The following inequalities hold: \begin{enumerate}
\item $\mu_1\geq \mu_2$,
\item $\mu_2\geq \mu_3$,
\item $\mu_3\geq \mu$.
\end{enumerate}
\end{LE}
\begin{proof}
{\bf (1)} If $\mu_1=\infty$, then the inequality $\mu_2\leq \mu_1$ holds clearly. Otherwise, $\mu_1$ is finite.

We claim that $\mu_1\geq 3$. For if not, then $\core(\mathcal{C})$ would have a cover $\{u,v\}$ of size $2$ that is monochromatic in some proper $2$-vertex-coloring of $G$. Clearly, $u,v$ must be from different connected components of $G$. Let $v'$ be an element in $V$ such that $\{v,v'\}$ is an edge of $G$. Then every member of $\core(\mathcal{C})$ containing $v'$ does not contain $v$ so it must contain $u$. Subsequently, by \Cref{setcore0}, $u$ and $v'$, and therefore $u$ and $v$, are from the same connected component of $G$, a contradiction. 

We prove by induction on $\mu_1$ that $\mu_2\leq \mu_1$. If $\mu_2\leq 3$, then the inequality follows from the inequality $\mu_1\geq 3$ that we just showed. For the induction step, we assume that $\mu_2> 3$. Suppose $U_1\cup \cdots\cup U_{\mu_1}$ is a cover of $\core(\mathcal{C})$, and let $\mathcal{C}':=\mathcal{C}\setminus U_{\mu_1}/V_{\mu_1}$, which is also a clean tangled clutter. As $\mu_2>3$, it follows that $G(\mathcal{C}')$ has the same connected components as $G$ except for $U_{\mu_1}\cup V_{\mu_1}$.

By the induction hypothesis, $\mu_2(\mathcal{C}')\leq \mu_1(\mathcal{C}')$. On the one hand, $\mu_1(\mathcal{C}')\leq \mu_1-1$ as $U_1\cup \cdots\cup U_{\mu_1-1}$ is a cover of $\core(\mathcal{C}')$; this is because $\core(\mathcal{C}')\subseteq \core(\mathcal{C})\setminus U_{\mu_1}/V_{\mu_1}$ by \Cref{deletion-contraction} and $U_1\cup \cdots\cup U_{\mu_1-1}$ is clearly a cover of the latter. On the other hand, $\mu_2(\mathcal{C}')\geq \mu_2-1$ as any monochromatic cover $U$ of $\mathcal{C}'$ yields a monochromatic cover of $\mathcal{C}$, namely $U\cup U_{\mu_1}$, intersecting only $1$ more component of $G$. Thus, $\mu_1-1\geq \mu_1(\mathcal{C}')\geq \mu_2(\mathcal{C}')\geq \mu_2-1$, implying in turn that $\mu_1\geq \mu_2$, thereby completing the induction step.

{\bf (2)} If $\mu_2=\infty$, then we are done. Otherwise, suppose $U_1\cup \cdots \cup U_{\mu_2}$ is a monochromatic cover of $\mathcal{C}$. Suppose for a contradiction that $U_1\cup \cdots \cup U_{\mu_2-1} \cup U_{\mu_2}$ is not irreducible, say $U_1\cup \cdots\cup U_{\mu_2-1}\cup V_{\mu_2}$ is also a cover of $\mathcal{C}$. Then $U_1\cup \cdots\cup U_{\mu_2-1}$ must be a cover of $\core(\mathcal{C})$, implying in turn that $\mu_2-1\geq \mu_1$, thus contradicting (1).

{\bf (3)} 
If $\mu_3=\infty$, we are done. Otherwise, the inequality follows from \Cref{irred-mono}.
\end{proof}

We are ready to prove the promised relation, that $\mu=\lambda$.

\begin{proof}[Proof of \Cref{mono-connectivity}]
By \Cref{mu-mu-mu}, $\mu_1\geq \mu_2\geq \mu_3\geq \mu$. We also know that $\mu\geq \mu_1=\lambda$, thus $\mu =\lambda$, as required.
\end{proof}

\section*{Acknowledgements}

This work is supported by ONR grant 14-22-1-2528 and EPSRC grant EP/X030989/1. We would like to thank an anonymous referee for a very helpful report on a precursor of this manuscript and suggesting the definitions of $\mu_1(\mathcal{C}),\mu_2(\mathcal{C})$, and the proof of \Cref{mu-mu-mu} part (1). 

\paragraph{Data Availability Statement.} No data are associated with this article. Data sharing is not applicable to this article.

{\small 
\bibliographystyle{alpha}
\bibliography{ideal-size-references}}
 
\end{document}